\documentclass[letterpaper, 10 pt, conference]{ieeeconf}  

\IEEEoverridecommandlockouts                              

\usepackage{graphics} 
\usepackage{epsfig} 
\usepackage{mathptmx} 
\usepackage{times} 
\usepackage{amsmath} 
\usepackage{amssymb}  
\usepackage{xcolor} 

\title{\LARGE \bf
Optimization of Hydrogen Blending in Natural Gas Networks for Carbon Emissions Reduction
}

\author{Mo Sodwatana$^{1}$, Saif R. Kazi$^2$, Kaarthik Sundar$^3$, and Anatoly Zlotnik$^2$
\thanks{*This study was funded by the U.S. Department of Energy's Advanced Grid Modeling (AGM) projects ``Joint Power System and Natural Gas Pipeline Optimal Expansion'' and ``Dynamical Modeling, Estimation, and Optimal Control of Electrical Grid-Natural Gas Transmission Systems'', as well as LANL Laboratory Directed R\&D Project ``Efficient Multi-scale Modeling of Clean Hydrogen Blending in Large Natural Gas Pipelines to Reduce Carbon Emissions''.  Research conducted at Los Alamos National Laboratory is done under the auspices of the National Nuclear Security Administration of the U.S. Department of Energy under Contract No. 89233218CNA000001.}
\thanks{$^{1}$Mo Sodwatana is a Ph.D. student in the Department of Energy Science and Engineering, Stanford University, Stanford, CA 94305, USA
        {\tt\small jarupas@stanford.edu}}%
\thanks{$^{2}$Saif R. Kazi and Anatoly Zlotnik are in the Applied Mathematics \& Plasma Physics Group, Los Alamos National Laboratory, Los Alamos, NM 87545, USA
        {\tt\small \{skazi,azlotnik\}@lanl.gov}}%
\thanks{$^3$Kaarthik Sundar is in the Information Systems \& Modeling Group, Los Alamos National Laboratory, Los Alamos, NM 87545, USA
        {\tt\small \{kaarthik\}@lanl.gov}}%
}

\begin{document}

\maketitle
\thispagestyle{empty}
\pagestyle{empty}

\begin{abstract}

We present an economic optimization problem for allocating the flow of natural gas and hydrogen blends through a large-scale transportation pipeline network. Physical flow of the gas mixture is modeled using a steady-state relation between pressure decrease and flow rate, which depends on mass concentration of the constituents as it varies by location in the network.  The objective reflects the economic value provided by the system, accounting for delivered energy in withdrawn flows, the cost of natural gas and hydrogen injections, and avoided carbon emissions.  The problem is solved subject to physical flow equations, nodal balance and mixing laws, and engineering inequality constraints. The desired energy delivery rate and minimum hydrogen concentration can be specified as upper and lower bound values, respectively, of inequality constraints, and we examine the sensitivity of the physical pressure and flow solution to these parameters for two test networks.  The results confirm that increasing hydrogen concentration requires greater energy expended for compression to deliver the same energy content, and the formulation could be used for valuation of the resulting mitigation of carbon emissions.  


\end{abstract}

\section{INTRODUCTION}

The United States bulk power grid increasingly relies on natural gas (NG) fueled generation and renewable energy, while planning for reduced dependence on fossil fuels.  The blending of hydrogen (H$_2$) generated using clean energy into natural gas delivery systems is proposed to support this energy transition while using the capital investments in existing pipelines for their entire lifespan  \cite{witkowski2018analysis}. The use of NG-H$_2$ mixtures has the potential to reduce greenhouse gas emissions from power generation as well as residential end-use \cite{raju2022}.   The ability to inject hydrogen produced using renewable energy into gas pipelines will also provide operational flexibility and storage capacity for the power grid, and takes advantage of an existing infrastructure that would otherwise be stranded with increasing electrification \cite{anderson2004harvesting}. There are however system integration issues and associated costs of such hydrogen utilization that affect pipeline transport, which must be considered using nontrivial extensions of methods developed for homogeneous gas transport \cite{chaczykowski2018gas}. 

Blends of hydrogen and natural gas exhibit different physical flow properties than pure hydrogen or natural gas. Several studies have examined the effects of the properties of the blended gas, such as density, viscosity, phase interactions, and energy densities, on the pipeline network and end-use applications \cite{mahajan2022,melaina2013blending}. In particular, when blends of hydrogen and natural gas are injected into a pipeline, without changes to operating setpoints, there is a reduction in transported gas pressure and energy delivered downstream. Whereas hydrogen has a higher calorific value by mass than methane, it has a lower density at given pressures. In a recent case study, the energy quantity transported for the same pressure ratio was reduced by 4\% and 14\% for a 10\% and a 40\% hydrogen blend, respectively, assuming equipment compatibility \cite{bainier2019impacts}. Moreover, due to the pressure reduction, the energy required for compression increases by 7\% and 30\% for the respective blend percentages.  Other recent technical studies focus on the feasibility and safety of injection and pressurization \cite{schuster2020centrifugal}. 

For natural gas transport, canonical problems utilize steady-state optimization to evaluate capacity \cite{koch2015evaluating} and determine economically optimal allocation \cite{rudkevich2017hicss}.   In these problems, the gas is assumed to be chemically homogeneous.  However, composition tracking may be important when a pipeline has multiple receipt points for species with different calorific values \cite{hante2019complementarity}.  Assuming that technical issues related to sealing, compression, and end-use appliances can be resolved and that engineering limitations on quantities such as minimum and maximum pressures and compressor energy can be specified \cite{erickson2022importance}, the design, operation, and economics of H$_2$-NG blend pipelines must be considered.  Energy from hydrogen gas produced by electrolysis is much more costly today than prevailing natural gas prices \cite{lee2021integrative}.  Quantifying the cost of avoided carbon emissions resulting from H$_2$ blending into NG pipelines, as well as the resulting change in energy delivery capacity, is therefore of particular interest.

In this study, we formulate an economic optimization problem for determining a feasible flow allocation that maximizes economic benefit for users of a pipeline system that transports a blend of two significantly different gases.  The formulation is effectively a single auction market mechanism in which suppliers offer natural gas or hydrogen at a given price per mass flow rate, and consumers bid for deliveries at a given price per energy content.  In addition, consumers can individually provide a bid for carbon mitigation at a given price per mass of carbon dioxide emissions avoided at that location.  This gives a specific value to the amount of energy delivered in the form of hydrogen that otherwise would arrive in the form of natural gas. Control variables in the formulation include optimized injections of pure NG or H$_2$, withdrawal rates of the mixture, and compressor operating setpoints.  The optimization is solved subject to physical flow equations, nodal balance and mixing laws, and engineering inequality constraints.

The rest of the manuscript is structured as follows.  In Section \ref{sec:network_model}, we specify network modeling and physical flow equations for a pipeline that transports a spatially inhomogeneous mixture of two gases arising from distributed injections of these constituents. In Section \ref{sec:optimization}, we define an objective function and additional inequality constraints that include engineering limitations, and formulate the optimization problem.  We present the results of sensitivity analyses performed for two case studies in Section \ref{sec:computational}, discuss the implications and potential follow-on studies in Section \ref{sec:discussion}, and briefly conclude in Section \ref{sec:conclusion}.

\section{Heterogeneous Gas Pipeline Network Model} \label{sec:network_model}

We consider a gas pipeline network that is represented using a directed graph with junctions \(j\in\mathcal{V}\) that are connected by pipes \((i,j)\in\mathcal{E}\).  Another set $\mathcal{C}$ of node-connecting elements is the collection of compressors, which are used to boost gas pressure.  Following previously developed notation \cite{rudkevich2017hicss}, we also specify a set \(\mathcal{G}\) of gNodes, where each gNode \(m \in \mathcal{G}\) represents a user of the pipeline that is associated to a physical node \(j(m) \in \mathcal{V}\). Note that more than one gNode can be located at a physical node.  Physical nodes with associated gNodes can be either supply or withdrawal nodes, because there cannot, for example, be one gNode that injects hydrogen while another makes a withdrawal at the same physical node. The sets of components that constitute the network are denoted by

\vspace{1ex}

\hspace{-3ex} \begin{tabular}{ll}
\( j \in \mathcal{V}\)  &\:  set of all physical nodes,\\
\( (i,j) \in \mathcal{E}\)  &\: set of edges representing pipes,\\
\( (i,j) \in \mathcal{C}\)  &\: set of edges with compressors,\\
\( j \in \mathcal{V}_s\)  &\:  set of slack physical nodes, subset of \(\mathcal{V}\),\\
\( m \in \mathcal{G}_s^{H_2}\)  &\:  set of gNodes that inject hydrogen,\\
\( m \in \mathcal{G}_s^{NG}\)  &\: set of gNodes that inject natural gas,\\
\( m \in \mathcal{G}_d\)  &\: set of gNodes that withdraw gas.
\end{tabular}
\vspace{1ex}

\noindent The physical state of the network is defined by the gas pressure $P_j$ and mass fraction of hydrogen $\gamma_j$  for each junction $j \in \mathcal{V}$, and the total mass flow $\phi_{ij}$ and mass fraction of hydrogen  $\gamma_{ij}$ on each pipe $(i,j)\in\mathcal{E}$.  Our key structural assumption is that flow directions are fixed \emph{a priori}. We suppose that the control variables available to the operator and users of the pipeline system are the compressor ratio $\alpha_{ij}$ of each compressor $(i,j)\in\mathcal{C}$, the hydrogen gas supply $s_m^{H_2}$ at gNode $m \in \mathcal{G}_s^{H_2}$, the natural gas supply $s_m^{NG}$ at gNode $m \in \mathcal{G}_s^{NG}$, and the withdrawal flow $d_m$ of the gas mixture from gNode $m \in \mathcal{G}_d$.  Diagrams of test networks used in our study are found in Figs. \ref{fig:singlepipe} and \ref{fig:8node}.  We now describe the relations between these quantities and modeling parameters, which represent the function of the network.


\subsection{Pipe Equations} \label{sec:pipe_equations}

We suppose that in steady-state flow, the hydrogen concentration is uniform along each edge $(i,j)$.  We use the Weymouth equation to model the relation between pressures at the endpoints of the pipe and the flow through it \cite{rios2015optimization}.  This equation is derived from the momentum conservation law in the Euler equations for one-dimensional turbulent flow in a pipe. The relation is given as

\begin{subequations} \label{eq:pipe_equations}
\begin{equation} \label{eq:weymouth}
    P_{i}^{2} - P_{j}^{2} = \frac{\lambda_{ij}L_{ij}}{D_{ij}A_{ij}^{2}} V_{ij} \phi_{ij} \left|\phi_{ij} \right|  \quad \forall (i,j) \in \mathcal{E}
\end{equation}

\noindent where \(V_{ij}\) [(m/s)$^2$] is the squared speed of sound in the blended gas,

\begin{equation} \label{eq:soundspeed}
    V_{ij} = \gamma_{ij}a_{H_{2}}^{2} + (1-\gamma_{ij})a_{NG}^{2} \quad \forall (i,j) \in \mathcal{E}.
\end{equation}
\end{subequations}

We approximate the state equation using the ideal gas law, so that partial pressures are additive.  The squared speed of sound in the blended gas is approximated as the linear combination of the squared wave speeds \(a_{NG}^2\) and \(a_{H_2}^2\) [m/s] in natural gas and hydrogen.  The specific gas constants of hydrogen and methane are 4.116 and approximately 0.478, respectively \cite{nesbitt2011handbook}, so that we expect the wave speed of hydrogen to be a factor of approximately 2.93 greater than that of natural gas.  Thus in this study, we use \(a_{NG} \approx 370 m/s\) and \(a_{H_2} \approx 1090 m/s\). Using the ideal gas approximation significantly simplifies our exposition of economic optimization. The assumptions used to obtain equations \eqref{eq:weymouth}-\eqref{eq:soundspeed} can in principle be relaxed to extend the results to the regime of non-ideal gases, which better approximates the conditions of gas transportation pipelines \cite{srinivasan2022numerical}.


\subsection{Nodal Equations} \label{sec:node_equations}

At every physical node \(j\), the net mass flow through the node and the net injection into the node must be balanced. We impose mass balance equations on natural gas and hydrogen, given by
\begin{subequations}
\begin{equation} \label{eq:ngflowbalance}
\begin{split}
    (1-\gamma_j) \sum_{k\in \partial_j^-} \phi_{jk} - \sum_{i\in \partial_j^+} (1-\gamma_{ij}) \phi_{ij} \qquad \qquad \\ \qquad \qquad = \sum_{m\in \partial_j^g} s_m^{NG} -  (1-\gamma_j)\sum_{m\in \partial_j^g} d_m,
\end{split}  \quad \forall j \in \mathcal{V},
\end{equation}
\begin{equation} \label{eq:h2flowbalance}
    \!\! \gamma_{j} \sum_{k\in\partial_{j}^{-}} \phi_{jk} - \sum_{i\in\partial_{j}^{+}} \gamma_{ij}\phi_{ij} = \! \sum_{m\in\partial_{j}^g} s_{m}^{H_2} - \gamma_{j} \! \sum_{m\in\partial_{j}^g} d_{m} \quad \forall j \in \mathcal{V}.
\end{equation}
\end{subequations}

\noindent Here \(\partial_{j}^{+}\) and \(\partial_{j}^{-}\) are the sets of nodes connected to node $j$ by incoming and outgoing edges, respectively.  Adding together the two equations \eqref{eq:ngflowbalance} and \eqref{eq:h2flowbalance} imposes the total nodal mass balance.  To enforce appropriate continuity in the concentration from a node to an outgoing edge, the following constraint is imposed at the node-to-edge interface:
\begin{equation}  \label{eq:continuity}
    \gamma_{i} = \gamma_{ij}  \quad \forall (i,j) \in \mathcal{E},
\end{equation}
where the hydrogen concentration at the edge \((i,j)\) leaving node \(i\) equals the concentration at node \(i\). We do not impose continuity from \((i,j)\) to the incoming node \(j\) as the concentration at \(j\) depends on that of all incoming edges. Finally, we suppose that at each slack node \(j \in \mathcal{V}_s\), the pressure is maintained at a nominal value \(\sigma_j\):
\begin{equation} \label{eq:slack_pressure}
    P_{j} = \sigma_j \quad \forall j \in \mathcal{V}_s.
\end{equation}
The slack node is typically used to represent a large source of natural gas such as a processing plant or storage facility.  Standard boundary conditions for flow network simulation, including gas pipelines, require at least one slack node for well-posedness, and we use this convention here.


\subsection{Compressor Modeling} \label{sec:comp_modeling}

Gas transmission pipelines are constructed with compressor stations, which are complex facilities that may have multiple compressor machines.  For the purpose of large-scale system modeling, we suppose that the action of a compressor station $(i,j)\in\mathcal{C}$ is aggregated as a pressure boost ratio $\alpha_{ij}$, which acts as
\begin{equation} \label{eq:comp_boost}
    P_{j}^{2} = \alpha_{ij}^{2} P_{i}^{2}   \quad \forall (i,j) \in \mathcal{C}.
\end{equation}

We consider the power used for gas compression as an important factor in our study, because increasing hydrogen fraction increases the amount of compression work required to transport a given amount of energy in the form of the gas blend.  Following standard practice \cite{menon2005gas}, the power \(W_c\) used to drive the compressor is formulated as
\begin{equation} \label{eq:comp_power}
    W_c = \left( \frac{286.76\cdot (\kappa_{ij}-1)\cdot T}{G_{ij}{\kappa_{ij}}} \right) \left( \alpha_{ij}^{m}-1 \right) \left|\phi_{ij}\right|, \quad \forall (i,j) \in \mathcal{C}
\end{equation}
where \(m = (\kappa_{ij} - 1)/\kappa_{ij}\). Here we use \(\kappa_{ij}\) to denote the specific heat capacity ratio for the mixed gas and \(G_{ij}\) to denote the specific gravity ratio. We approximate these ratios for the blend using linear combinations of the specific ratios of each gas by
\begin{subequations}  \label{eq:gasratios}
\begin{align}
    \kappa_{ij} = \kappa_{H_{2}}\gamma_{ij} + \kappa_{NG}(1-\gamma_{ij}) \quad \forall (i,j) \in \mathcal{C}, \\
    G_{ij} = G_{H_{2}}\gamma_{ij} + G_{NG}(1-\gamma_{ij}) \quad \forall (i,j) \in \mathcal{C}.
\end{align}
\end{subequations}

\noindent In this study, we use \(\kappa_{NG} = 1.304\), \(\kappa_{H_{2}} = 1.405\), \(G_{NG} = 0.5537\), and \(G_{H_{2}} = 0.0696\). \(T\) [K] is the compressor suction temperature at which gas is being compressed, and is kept at 288.7 K in our simulations.  Here \(\phi_{ij}\) [kg/s] is the flow rate along the compressor edge \((i,j)\).


\subsection{Carbon Emissions Offset} \label{sec:emissions_modeling}

A key innovation of our study is to include carbon emissions mitigation based on the value of carbon displacement that is submitted as part of a bid by each consumer of energy to the optimization-based auction market mechanism.  The emissions offset \(E_m\) of a specific consumer $m\in\mathcal{G}_d$ is formulated as
\begin{equation}  \label{eq:carbon_offset}
    E_m = d_m \gamma_{j(m)} \cdot \frac{R_{H_2}}{R_{NG}}   \cdot \zeta_{NG},
\end{equation}
where \(\zeta_{NG} \) is the ratio of molecular weights of carbon dioxide and natural gas, and is approximately 44/18.  Here the constants \(R_{NG}\) and \(R_{H_2}\) are the calorific values for natural gas and hydrogen, respectively, for which we use values of \(R_{H_{2}} = 141.8 MJ/kg\) and \(R_{NG} = 44.2 MJ/kg\). The emissions term \eqref{eq:carbon_offset} denotes the amount of carbon dioxide [kg/s] that was \emph{not} emitted because hydrogen was burned instead of natural gas to produce a given amount of energy using the delivered flow \(d_m\) [kg/s] at concentration \(\gamma_{j(m)}\). Here, we examine the emissions from the energy replaced, rather than the mass replaced, from burning hydrogen, and therefore we include the ratio of calorific value of hydrogen to natural gas.






\section{Optimization Formulation} \label{sec:optimization}
We now describe the economic optimization formulation for constrained heterogeneous gas transport.  In addition to the physical and network modeling in Section \ref{sec:network_model}, we define an objective function, inequality constraints that specify engineering limitations, and parameters that constitute an economic bid structure.


\subsection{Economic Value Objective Function} \label{sec:objective}

Our objective is to maximize the economic value of transporting gas between suppliers of constituent gases as commodities and consumers who purchase energy and carbon emissions offsets. The economic value produced by the pipeline is the sum over all gNodes \(m\) of payments by consumers for delivered blended gas minus purchases from suppliers of received gas constituents. Suppliers place offer prices for natural gas, \(c_m^{NG}\) [\$/kg], and hydrogen, \(c_m^{H_2}\) [\$/kg], at the supplying gNodes, while off-takers place bids for the energy content of the blended gas, \(c_m^d\) [\$/MJ], at the off-taking gNodes. Off-takers can also place a value on the carbon emissions avoided, \(c_m^{CO_2}\) [\$/kg]. The offer price is in terms of mass while the bid price is in energy units, which represents how suppliers and off-takers consider the value of gas differently in a market that includes the cost of emissions. We also include the work done to compress gas by each compressor \(c\in\mathcal{C}\), which is denoted by \(W_c\) in the economic value formulation.  The economic value objective function is then expressed as
\begin{equation}\label{eq:obj}
\begin{split}
    J_{EV} = \sum_{m\in G} \bigg(c_m^{d} d_m (R_{H_2}\gamma_{j(m)} + R_{NG} (1-\gamma_{j(m)}))  \qquad \\    - c_m^{H_{2}s} s_m^{H_2}  - c_m^{NGs} s_m^{NG} +  c_m^{CO_2} E_m \bigg)  - \eta \sum_{c\in C} W_c
\end{split}
\end{equation}
where \(\eta\) [\$/kw-s] is a conversion factor used to define the economic cost of applied compressor power. \(R_{NG}\) and \(R_{H_2}\) are the calorific values of burning natural gas and hydrogen, respectively. The calorific value of the blended gas is the linear combination of the calorific values of the two gases with respect to \(\gamma_{j(m)}\), the mass fraction of hydrogen at \(j\). In our computational case studies, \(\eta = \)\$0.13/3600kw-s, and \(R_{H_{2}} = 141.8 MJ\)/kg and \(R_{NG} = 44.2 MJ\)/kg as specified in Section \ref{sec:emissions_modeling}. Here \(d_m\) [kg/s] is the mass flow rate of delivered blended gas while \(s_m^{NG}\) and \(s_m^{H_2}\) [kg/s] are the respective mass flow rates of natural gas and hydrogen at the supply gNodes.


\subsection{Pressure, Compressor, and Concentration Limits} \label{sec:engineering_limits}

We suppose that minimum pressure limits as well as minimum and maximum hydrogen concentration limits may be specified at each node:
\begin{subequations}
\begin{align} 
    P_{j}^{min} \le P_{j}  \quad & \forall j \in \mathcal{V}, \label{eq:minpressure} \\
    \gamma_{j}^{min} \le \gamma_{j} \le \gamma_{j}^{max} \quad & \forall j \in \mathcal{V}. \label{eq:conclimits}
\end{align}
\end{subequations}
In addition, the discharge node of each compressor has a maximum allowable operating pressure, and each compressor has a maximum boost ratio:
\begin{subequations}
\begin{align} 
    \alpha_{ij}P_{i} \le P_{ij}^{max}  \quad & \forall (i,j) \in \mathcal{C}, \label{eq:maxpressure} \\
    1 \le \alpha_{ij} \le \alpha^{max}_{ij} \quad & \forall (i,j) \in \mathcal{C}. \label{eq:complimits}
\end{align}
\end{subequations}

\subsection{Supply and Demand Limits} \label{sec:economic_limits}

The supplies of natural gas and hydrogen at gNodes \(m\) are positive, and are constrained by upper bounds of the form
\begin{subequations} \label{eq:supplylimits}
\begin{align}
    0\le s_{m}^{NG} \le s_{m}^{max,NG} \quad & \forall m \in \mathcal{G}_s^{NG}, \\
    0\le s_{m}^{H_2} \le s_{m}^{max,H_2} \quad & \forall m \in \mathcal{G}_s^{H_2}.
\end{align}
\end{subequations}


\noindent We suppose that outflows $d_m$ from the network are positive. The upper bound of the demand is in terms of energy, with calorific values as conversion from mass flow to energy flow.  Energy deliveries are either optimized or fixed, corresponding to gNode sets $\mathcal{G}_{d,o}$ and $\mathcal{G}_{d,f}$,  with constraints of form
\vspace{-2ex}
\begin{subequations}\label{eq:energydemand}
\begin{align}
   &\!\!\!\!\! 0\le d_{m} \left( R_{H_{2}}  \gamma_{j(m)} + R_{NG} (1\!-\gamma_{j(m)}) \right) \le g_{m}^{max}, \,\,   \forall m \!\in \mathcal{G}_{d,o}, \!\!\! \label{eq:demand_opt} \\
  &0\le d_{m} \left( R_{H_{2}}  \gamma_{j(m)} + R_{NG} (1-\gamma_{j(m)}) \right) = \bar{g}_{m}, \,\,   \forall m \in  \mathcal{G}_{d,f}. \!\!\! \label{eq:demand_fixed}
\end{align}
\end{subequations}
Combining equations \eqref{eq:weymouth}-\eqref{eq:energydemand} yields the economic heterogeneous gas transport optimization problem:
\begin{equation} \label{prob:hgto}
\begin{array}{ll}
\!\!\!\! \mathrm{max}  & J_{EV} \triangleq \text{max economic value objective } \eqref{eq:obj} \\
\!\!\!\! \text{s.t.} & \text{pipe flow equation } \eqref{eq:pipe_equations} \\
& \text{NG nodal flow balance } \eqref{eq:ngflowbalance} \\
& \text{H\textsubscript{2} nodal flow balance } \eqref{eq:h2flowbalance} \\
& \text{concentration continuity } \eqref{eq:continuity} \\
& \text{slack pressure } \eqref{eq:slack_pressure} \\
& \text{compression } \eqref{eq:comp_boost} \\
& \text{pressure limits } \eqref{eq:minpressure},\eqref{eq:maxpressure} \\
& \text{H\textsubscript{2} concentration limits } \eqref{eq:conclimits} \\
& \text{compressor boost limits } \eqref{eq:complimits} \\
& \text{supply limits } \eqref{eq:supplylimits} \\
& \text{demand limits } \eqref{eq:demand_opt} \\
& \text{fixed consumption } \eqref{eq:demand_fixed} 
\end{array}
\end{equation}

\noindent We interpret the parameters $(g_{m}^{max},c_m^{d},c_m^{CO_2},\gamma_{j(m)}^{min})$ as a market bid by gas consumers, with a quantity (MJ) and price (\$/MJ) of delivered energy, a carbon offset price (\$/Kg), and minimum hydrogen mass fraction (\%). The decision variables in problem \eqref{prob:hgto} are given with units in Table \ref{tab:primal}.

\begin{table}[!t]
\vspace{2ex}
\centering \normalsize
\begin{tabular}{|l|c|c|c|c|}
  \hline
  Variables & Set & Units (SI) \\ 
  \hline
  $s_m^{H_2}$ & $m\in \mathcal{G}_s^{H_2}$ & kg$\cdot$ s$^{-1}$ \\
  $s_m^{NG}$  & $m\in \mathcal{G}_s^{NG}$ & kg$\cdot$ s$^{-1}$ \\
  $d_m$ & $m\in \mathcal{G}_d$ & kg$\cdot$ s$^{-1}$ \\
  $\alpha_{ij}$ & $(i,j)\in \mathcal{C}$  & - \\
  $\phi_{ij}$ & $(i,j)\in \mathcal{E}$  & kg$\cdot$ s$^{-1}$ \\
  $\gamma_{ij}$ & $(i,j)\in \mathcal{E} $ & - \\
  $\gamma_{j}$ & $j\in \mathcal{V} $ & - \\
  $P_{j}$ & $j\in \mathcal{V}$  & kg$\cdot$ m$^{-1}$$\cdot$ s$^{-2}$ \,\, (Pa) \\
  \hline
\end{tabular}
\caption{Collection of optimization variables for Problem \eqref{prob:hgto}. \label{tab:primal}}
\end{table}

\section{Case Studies} \label{sec:computational}

\begin{figure}[!t] 
    \vspace{-5ex}
    \includegraphics[width=\linewidth]{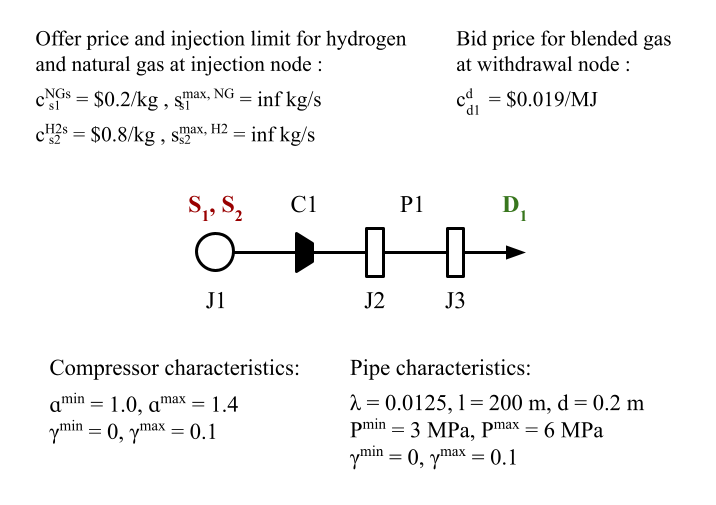}
    \centering
    \vspace{-6ex}
    \caption{Schematic of a single pipe network with compressor C1, natural gas and hydrogen injection at gNodes, S1 and S2, respectively, and a withdrawal gNode D1. The pipe and compressor characteristics are shown, as well as offer and bid prices at the supply and demand node, respectively.}
    \label{fig:singlepipe}
\end{figure}

We examine solutions to the problem \eqref{prob:hgto} using two test networks -- a single pipe with a compressor at the start (Fig. \ref{fig:singlepipe}) and an 8-node network with a loop and three compressors (Fig. \ref{fig:8node}). Problem \ref{prob:hgto} is solved in the Julia programming language v1.7.2 using the JuMP package v1.1.0 \cite{dunning2017jump}, which is an optimization modeling toolkit. We use IPOPT v1.0.2, a large scale optimization package \cite{wachter2006implementation}, as the nonlinear program solver. The case studies are computed on an AMD EPYC 7742 64-core processor with 16 GB of RAM. The solve times for both networks are mere seconds, and sensitivity analyses involving batches of optimization problems require approximately five minutes or less, depending on network size and other factors.

\subsection{Non-Dimensionalization and Rescaling}

We non-dimensionalize the governing equations prior to solving problem \eqref{prob:hgto} in order to avoid numerical issues \cite{srinivasan2022numerical}. In addition, we re-scale \eqref{eq:weymouth} because the wave speed \(V\) in blended gas  is orders of magnitude larger than other variables in the equation. Let \(\bar{P} = P / P_0\), \(\bar{L} = L/ l_0\), \(\bar{D} = D/ l_0\), \(\bar{A} = A/ A_0\), and \(\bar{\phi} = \phi/ \phi_0 = \phi / (\rho_0 u_0 A_0)\). Equation \eqref{eq:weymouth} then becomes
\begin{subequations} \label{eq:pipe_equations_dimensionless}
\begin{align}
    \bar{P}_{i}^{2} - \bar{P}_{j}^{2} = \frac{\lambda_{ij}\bar{L}_{ij}}{\bar{D}_{ij}\bar{A}_{ij}^{2}} \bar{V}_{ij}  \bar{\phi_{ij}} \left|\bar{\phi}_{ij} \right| \cdot \frac{u_0^2}{a_0^2} \cdot  \quad \forall (i,j) \in \mathcal{E}, \label{eq:weymouth_dimensionless} \\ 
    \bar{V}_{ij} \triangleq \frac{V_{ij}(\gamma_{ij})}{a_0^2} \quad \forall (i,j) \in \mathcal{E}. \label{eq:soundspeed_dimensionless}
\end{align}
\end{subequations}

\noindent We re-scale wave speed with a factor of \(a_0 = 635.06\) m/s. For the single pipe and 8-node cases, the nominal pressures are \(P_0 = 5 \) MPa and \(P_0 = 3.04\) MPa, respectively. The nominal length, area, density and velocity for both cases are \(l_0 = 5000\) m, \(A_0 = 1\) m$^2$, \(\rho_0 = P_0 / a_0^2\), and \(u_0 = \lceil a_0 \rceil / 300\), where \(a_0\) is the geometric mean of wave speeds, used in the re-scaling factor. We compute \(a_0\) as \(a_0 = \sqrt{a_{NG} \cdot a_{H_2}} \) where \(a_{NG} = 370\) m/s and \(a_{H_2} = 1090\) m/s are the wave speeds of NG and H$_2$, obtained by \(a_{NG} = \sqrt{RT/M_{NG}}\) and \(a_{H_2} = \sqrt{RT/M_{H_2}}\), respectively. Here, \(R = 8.314\) J/mol/K is the universal gas constant, and \(M_{NG} = 0.01737\) kg/mol and \(M_{H_2} = 0.002016\) kg/mol are molecular masses of NG and H$_2$.

\subsection{Single Pipe} \label{sec:singlepipe}

\vspace{-0.5ex}
For the single pipe test network, we examine how the solution changes due to 1) variation of the constraint bound $g_{m}^{max}$ for maximum energy demand \eqref{eq:demand_opt}; 2) a similar analysis with fixed nodal values of carbon emissions offsets $c_m^{CO_2} E_m$ added to the objective; 3) variation of the constraint bound $\gamma_{j}^{min}$ for minimum hydrogen concentration \eqref{eq:conclimits}; and 4) variation of the carbon emissions offset price $c_m^{CO_2}$. The structure and characteristics of the single pipe network are described in Figure \ref{fig:singlepipe}, and the four sensitivity analyses and their results are described in Sections \ref{sec:analysis1} through \ref{sec:analysis4}.

\subsubsection{Sensitivity analysis with respect to the maximum energy demand at the withdrawal gNode, \(g^{\max}_{D1}\), without carbon emissions offset values} \label{sec:analysis1}

\begin{figure}[!t]
    \includegraphics[width=\linewidth]{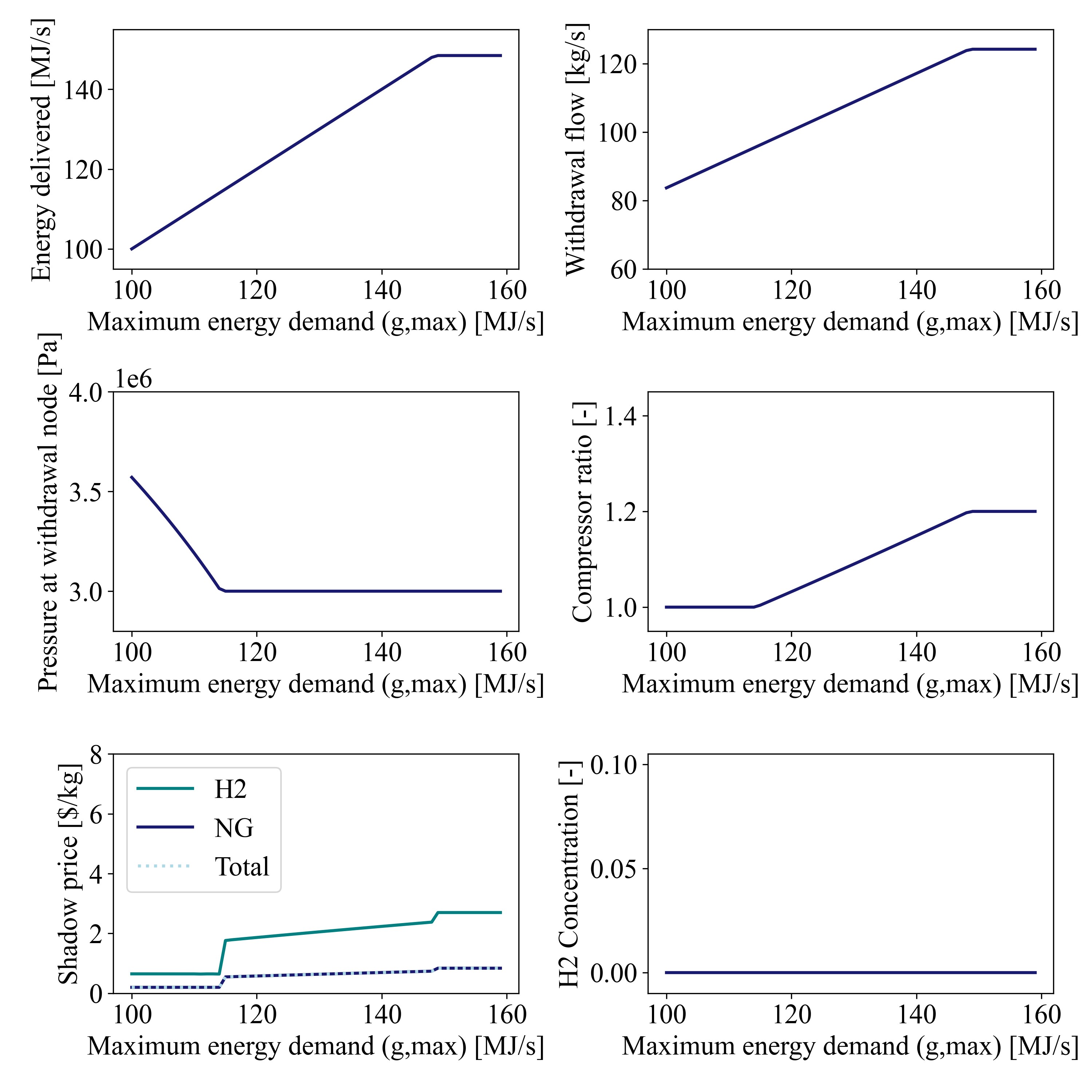}
    \centering
    \vspace{-4ex}
    \caption{Results of sensitivity analysis for the single pipe network in Section \ref{sec:analysis1} are shown as the change in physical variables at the withdrawal gNode D1, the shadow price at the physical node J3, and the compressor ratio at C1.  The maximum energy demand \(g_{D1}^{max}\) is varied from 100 to 160 MJ/s. There is no carbon emissions offset at the withdrawal node \((c_{D1}^{CO_2} \equiv 0)\).}
    \label{fig:singlepipe-energy-nocc}
\end{figure}

Here we vary the constraint bound \(g^{\max}_{D1}\) from 100 MJ/s to 160 MJ/s in increments of 1 MJ/s, and solve problem \eqref{prob:hgto} for each instance.  The variations in the solutions are shown in Figure \ref{fig:singlepipe-energy-nocc}, and can be divided into three regions where the transitions arise from activation of new binding constraints.  In the first region, where $100\leq g^{\max}_{D1} \leq 115$ MJ/s, the energy demand constraint binds because the system is able to meet the demand. The pressure at each node is within the defined bounds and there is no work done by the compressor. At 115 MJ/s, the pressure at the receiving node D1 hits the lower bound of 3 MPa, and therefore the compressor must operate. In the region where $ 115 \leq g^{\max}_{D1} \leq 150$ MJ/s, the compressor ratio increases to meet the rising energy demand. In the final region, where $150 \leq g^{\max}_{D1}$ MJ/s, energy and mass flow delivered taper off to become constant regardless of further increases in energy demand, with the compressor ratio binding at its limit of 1.4.

The shadow prices of H$_2$ and NG are shown to change with \(g^{\max}_{D1}\) as well in Figure \ref{fig:singlepipe-energy-nocc}.  These quantities are computed by the solver, and are given as the duals of the respective mass balance constraints at each physical withdrawal node. The shadow price of the mixture is the linear combination of the two shadow prices weighted by mass fraction ratios.  When a constraint becomes active and a phase transition occurs, the shadow prices adjust accordingly. As more power must be used for gas compression to deliver energy, there is a steady increase in the shadow prices. When the energy demand can no longer be met, there is a jump in the shadow prices.

\subsubsection{Sensitivity analysis with respect to the maximum energy demand at the withdrawal gNode, \(g^{\max}_{D1}\), with carbon emissions offset values} \label{sec:analysis2}

\begin{figure}[!t]
    \includegraphics[width=\linewidth]{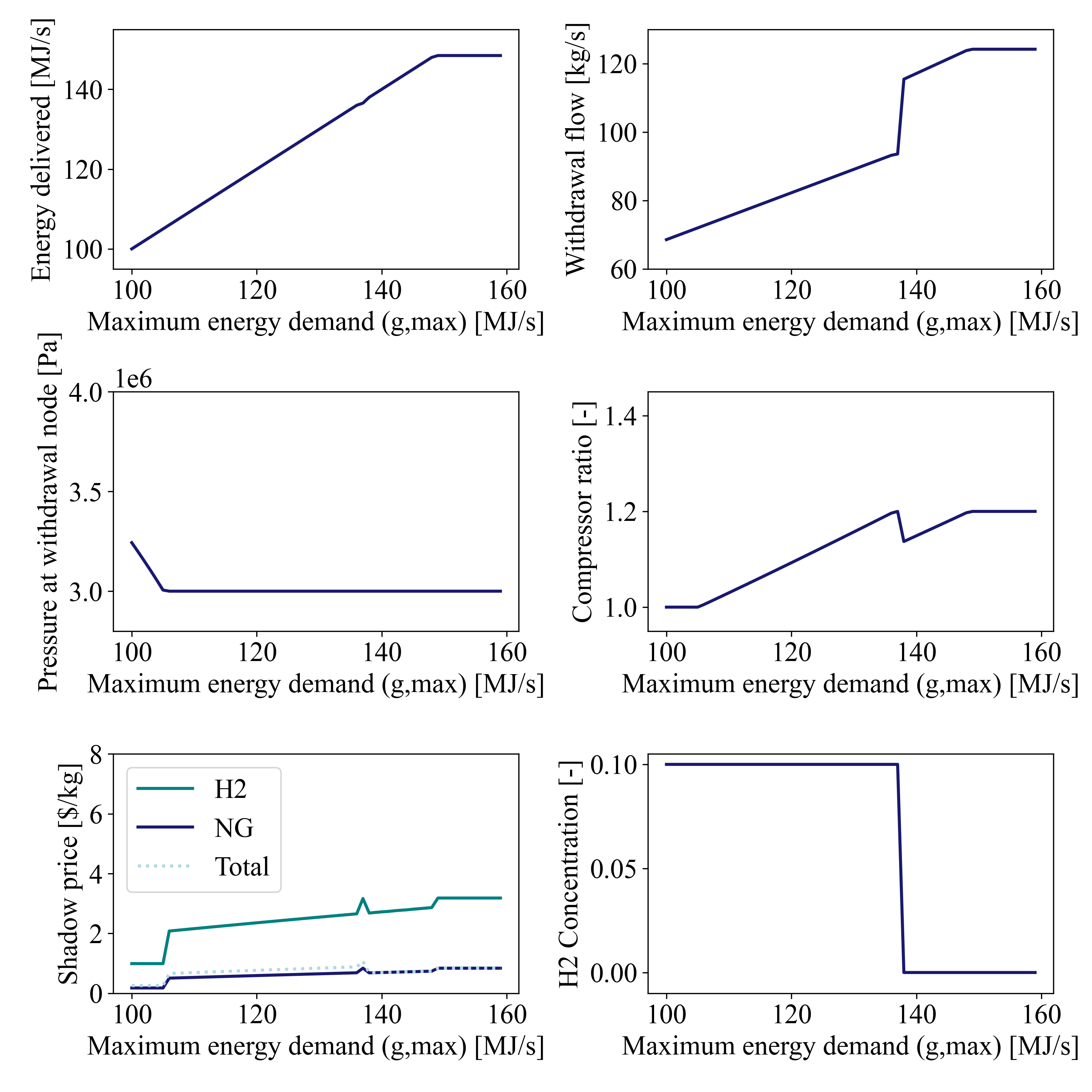}
    \centering
    \vspace{-4ex}
    \caption{Results of sensitivity analysis for the single pipe network in Section \ref{sec:analysis2} are shown as the change in physical variables at the withdrawal gNode D1, the shadow price at the physical node J3, and the compressor ratio at C1. The maximum energy demand \(g_{D1}^{max}\) is varied from 100 to 160 MJ/s when \(c_{D1}^{CO_2}=\$0.055/kg\) at the withdrawal node.}
    \label{fig:singlepipe-energy-wcc}
\end{figure}

Here we repeat the analysis in Section \ref{sec:analysis1}, and include a carbon emissions offset value where \(c_{D1}^{CO_2} = \$0.055\)/kg CO$_2$ avoided. The changes in the solution in this case are shown in Figure \ref{fig:singlepipe-energy-wcc}. The bound value \(g^{\max}_{D1}\) increases from 100 MJ/s to 160 MJ/s in increments of 1 MJ/s as with the previous case. The solutions can be categorized into four regions.  Providing a value of reducing carbon emissions incentivizes end-use of hydrogen, so the solution initially sees hydrogen injected at the maximum allowable concentration of 0.1. In this first region of four, where \(g^{\max}_{D1}\leq 105\) MJ/s, the maximum energy demand is the binding constraint. A transition occurs at 105 MJ/s where the pressure at junction J3 binds at the lower limit of 3 MPa. This pressure constraint is binding at a lower maximum energy demand compared to the previous case, because hydrogen blending reduces the pipeline pressure. For  $105 \leq g^{\max}_{D1} \leq  137$ MJ/s, there is compressor work done to ensure that the maximum energy delivery constraint is being met while pressure does not fall below the lower limit. In this region, the slope of the withdrawal flow rate becomes steeper because hydrogen in the gas mixture requires more natural gas by mass to substitute for the same energy. In region three, where  $137 \leq g^{\max}_{D1} \leq  150$ MJ/s, the solution switches to injecting pure NG to maintain allowable pressure while meeting increasing energy demands. This behavior is contrary to the previous case where the energy delivered tapers off.  In the final region, for $g^{\max}_{D1} \geq  150$  MJ/s, the energy delivered and withdrawal flow become constant as the system becomes congested cannot supply energy at the rate \(g^{\max}_{D1}\).

\begin{figure}[!t]
    \includegraphics[width=\linewidth]{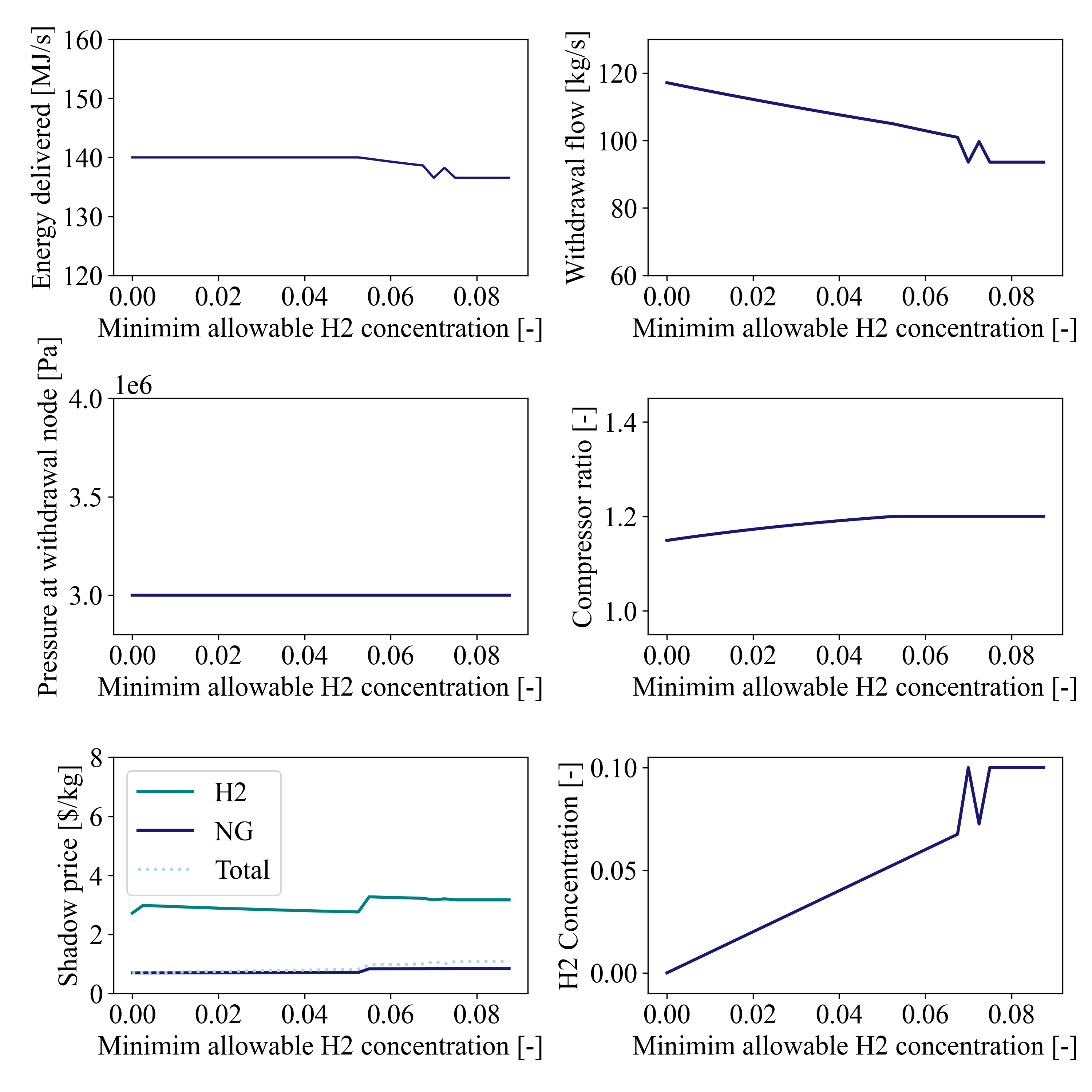}
    \centering
    \vspace{-4ex}
    \caption{Results of sensitivity analysis in Section \ref{sec:analysis3} with respect to the minimum allowable hydrogen concentration at the withdrawal gNode, D1. Here, the maximum energy demand is constant at 140 MJ/s.}
    \label{fig:singlepipe-conc-wcc}
\end{figure}

\subsubsection{Sensitivity analysis with respect to the minimum allowable H$_2$ concentration, \(\gamma_{J3}^{\min}\), at the withdrawal gNode with a carbon emissions offset} \label{sec:analysis3}

We examine the sensitivity of the solution to the minimum hydrogen concentration \(\gamma_{J3}^{\min}\) at the physical withdrawal node J3 corresponding to the withdrawal gNode D1. This scenario reflects instances where a consumer is must mitigate their carbon emissions by bidding for a minimum H$_2$ concentration at the delivery node. We increase \(\gamma_{J3}^{\min}\) from 0.0 to 0.1 at an increment of 0.025 while \(g_{D1}^{\max}\) is constant at 140 MJ/s.  Figure \ref{fig:singlepipe-conc-wcc} provides a summary of the results.  The solutions can be categorized in three regions.  For $0\leq\gamma_{J3}^{\min}\leq 0.053$, the energy demand  \(g_{D1}^{\max}\) can be met, and the actual delivered concentration increases with the minimum bound value.  The compressor ratio must also increase to keep the line pressure high enough to meet the minimum pressure bound at J3.  For $0.053\leq\gamma_{J3}^{\min}\leq 0.068$, The maximum pressure at J2 binds so that the pipe is congested, and increasing the minimum hydrogen concentration results in less energy delivered.  At $\gamma_{J3}^{\min}\approx 0.068$, the concentration binds at the maximum value.  The change in primal and dual solutions is non-monotone as concentration parameters change, which is consistent with prior observations on transport of heterogeneous gas blends.

\subsubsection{Sensitivity analysis with respect to the carbon emissions offset price $c_m^{CO_2}$ at the withdrawal gNode, D1} \label{sec:analysis4}

\begin{figure}[!t]
    \includegraphics[width=\linewidth]{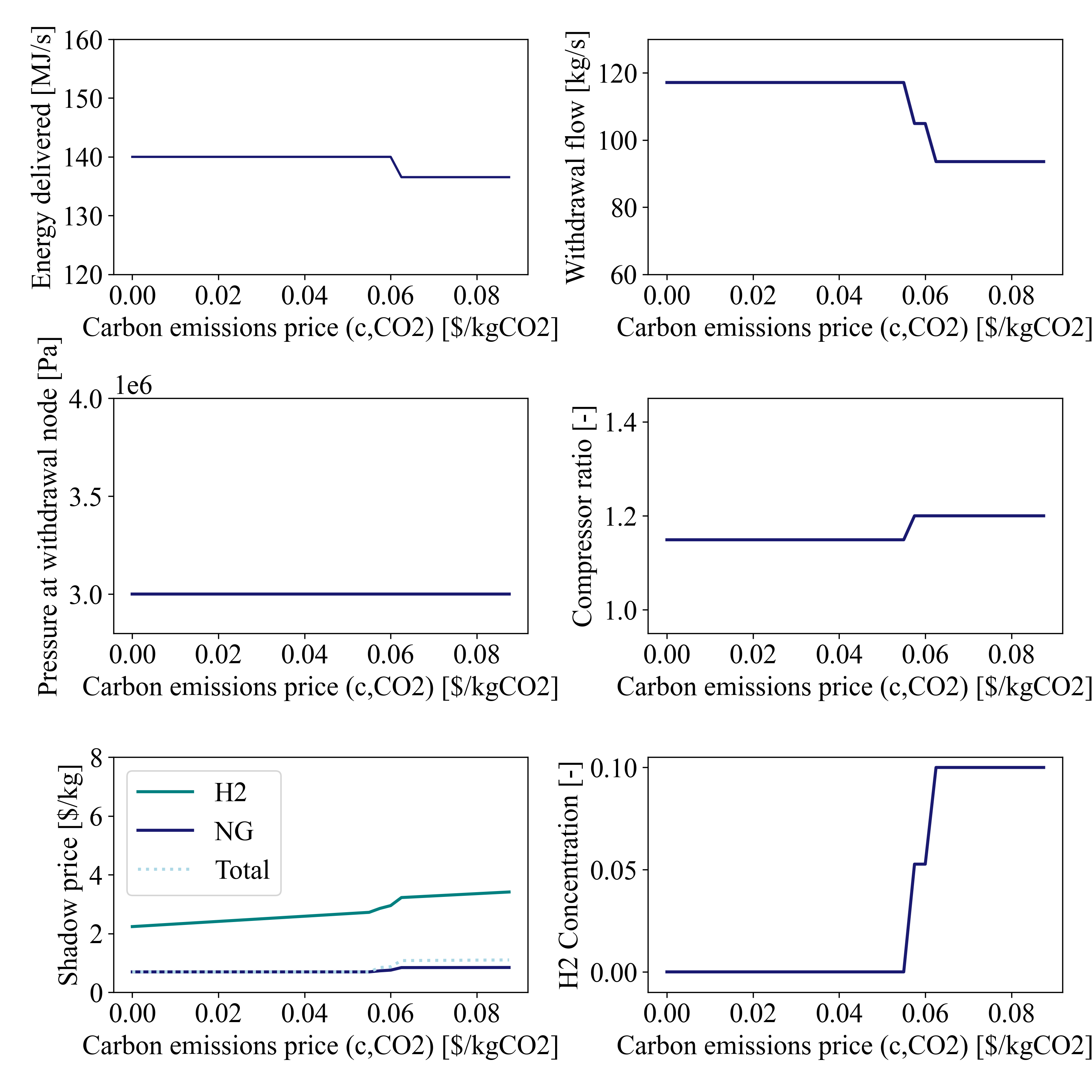}
    \centering
    \vspace{-4ex}
    \caption{Results of sensitivity analysis in Section \ref{sec:analysis4} with respect to the carbon emissions offset price $c_m^{CO_2}$ at the withdrawal gNode, D1. Here, the maximum energy demand is constant at 140 MJ/s.}
    \label{fig:singlepipe-carb-price}
\end{figure}

For this analysis, we vary the carbon emissions offset price $c_m^{CO_2}$ at the withdrawal gNode D1 from 0 to \$0.09 \$/kgCO$_2$.  The solutions can be categorized into three regions of interest.  The energy demand \(g_{D1}^{\max}\) is set to 140 MJ/s, and the minimum concentration \(\gamma_{J3}^{\min}\) is set to 0.  Figure \ref{fig:singlepipe-carb-price} shows a summary of the results.  When $c_m^{CO_2}\leq 0.055$, no H$_2$ is injected into the system, and $\gamma_{J3}\equiv 0$.  The energy delivered is binding at the maximum, and the flow is pure NG.  Though the shadow price of H$_2$ gradually increases to reflect its increased value, it is still not sufficient for injection of H$_2$ to add overall economic value for the system.  In the region $0.055 \leq c_m^{CO_2}\leq 0.065$, several transition occur.  First, the value of carbon mitigation surpasses the value of delivered energy, so that H$_2$ starts to be injected into the system, and the maximum pressure at node J2 binds so that the pipeline becomes constrained.  As a result, less energy is delivered than the request at gNode D1, the withdrawal flow drops, and more hydrogen is injected into the system since its carbon offset value continues to be greater than the value of delivered energy at the offer price.  Then in the region $c_m^{CO_2}\geq 0.065$, the maximum H$_2$ concentration constraint at the physical node J3 binds, and the only subsequent change is an increase in the shadow price of H$_2$ to reflect its increasing carbon mitigation value.

\subsection{8-Node Network} \label{sec:8node}

\begin{figure}[!t]
    \includegraphics[scale=0.35]{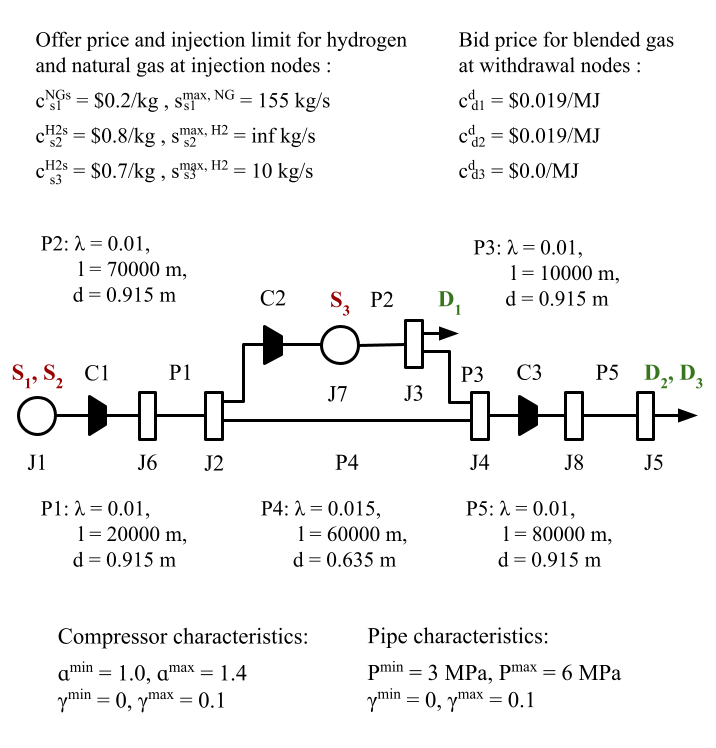}
    \centering
    \vspace{-6ex}
    \caption{Schematic of an 8-node network with one limited NG supplier at gNode S1, a supplier of unlimited H$_2$ at gNode S2, and a limited supplier of H$_2$ at S3. The limited supplier at S3 provides H$_2$ at a relatively cheaper price than the H$_2$ supplier at S2. At D1 and D2 \(c_{D1}^{CO_2}\) = \(c_{D2}^{CO_2}\) = \$0.055/kg CO$_2$ avoided. D3 is a fixed demand node with \(\bar{g}_m^{\max}\) = 100 MJ/s and \(c_{D3}^{CO_2}\) = 0.  gNodes S1 and S2 share the same physical node J1. gNodes D2 and D3 share J5 as the physical node. Pipe length, diameter, friction factor varies.}
    \label{fig:8node}
\end{figure}

\begin{figure}[!t]
    \includegraphics[width=\linewidth]{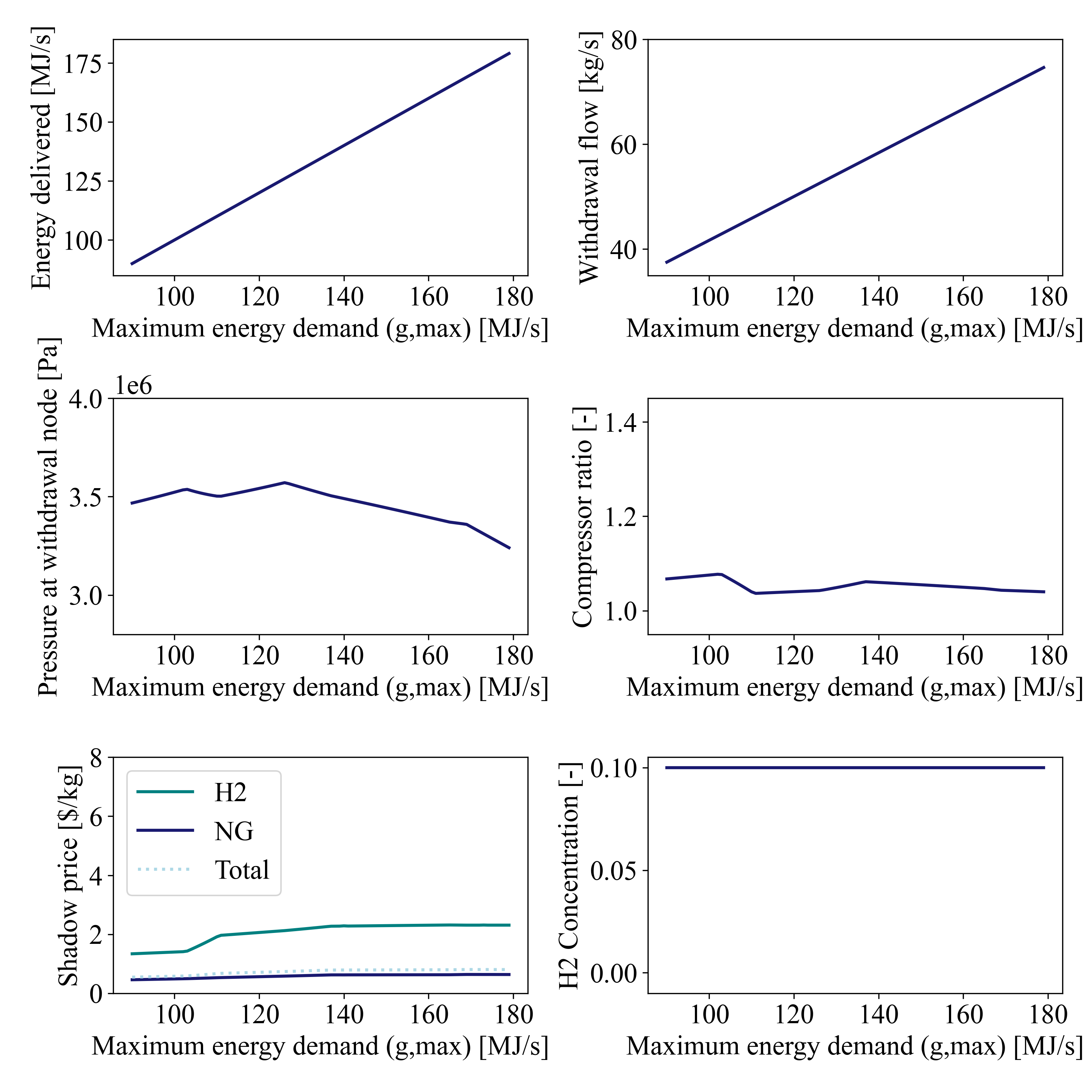}
    \centering
    \vspace{-4ex}
    \caption{Plots of the changing physical characteristics at the withdrawal gNode D1, the shadow price at the physical node J3, and the compressor ratio at C2 from a sensitivity analysis \(g_{D1,D2}^{max}\) on the 8-node pipe network. Here, \(c_{D1}^{CO_2}.\) = \$0.055/kg CO$_2$ avoided at the withdrawal gNode.}
    \label{fig:8nodepipe-energy-wcc-nodeD1}
    \vspace{-1ex}
\end{figure}

\begin{figure}[!t]
\vspace{-1ex}
    \includegraphics[width=\linewidth]{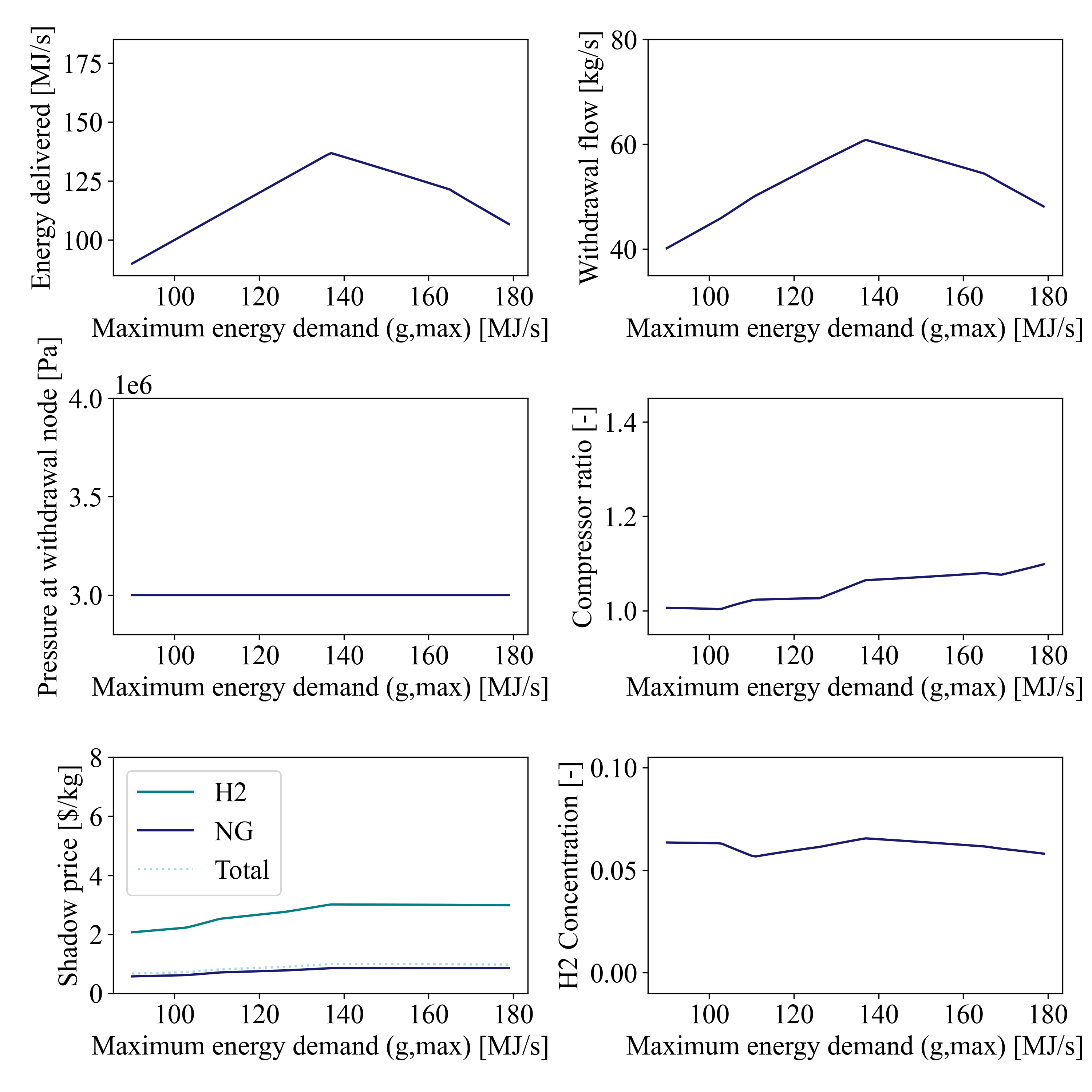}
    \centering
    \vspace{-4ex}
    \caption{Plots of the changing physical characteristics at the withdrawal gNode D2, the shadow price at the physical node J5, and the compressor ratio at C3 from a sensitivity analysis \(g_{D1,D2}^{max}\) on the 8-node pipe network. Here, \(c_{D2}^{CO_2}.\) = \$0.055/kg CO$_2$ avoided at the withdrawal gNode.}
    \label{fig:8nodepipe-energy-wcc-nodeD2}
\end{figure}

We use the 8-node network shown in Figure \ref{fig:8node} to demonstrate the behavior of solutions to optimization problem \eqref{prob:hgto} in the case of a more complex network topology.  First, we see that the computational implementation readily scales to looped topologies and multiple edges and junctions.  We also perform a sensitivity analysis by varying the maximum energy demands \(g_{D1}^{\max}\) and \(g_{D2}^{\max}\) at both withdrawal gNodes with a carbon emissions offset fixed at \(c_{D1}^{CO_2} = \$0.055\)/kg and the minimum H$_2$ concentrations fixed at \(\gamma_{J3}^{\min} = \gamma_{J5}^{\min} = 0.05\).  We increase \(g_{D1}^{\max}\) and \(g_{D2}^{\max}\) from 120 MJ/s to 180 MJ/s in increments of 1 MJ/s, and present the results seen at gNode D1 in Figure \ref{fig:8nodepipe-energy-wcc-nodeD1} and results for gNode D2 in Figure \ref{fig:8nodepipe-energy-wcc-nodeD2}. 
We observe two regions in the solution. In the first region where \(g_{D1}^{\max},g_{D2}^{\max} \le 140\) MJ/s, the solution can continually meet both \(g_{D1}^{\max}\) and \(g_{D2}^{\max}\). In the second region, \( g_{D1}^{\max},g_{D2}^{\max} \geq 140\) MJ/s, the energy and the flow delivered at D2 gradually declines because the system becomes congested due to a constraint on how much gas can be delivered to physical node J5. Because D2 shares a physical node with D3, a gas consumer with a fixed demand, the solution has to reduce supply to D2. Meanwhile, D1 receives uninterrupted supply of the blended gas because of its proximity to the injection site. 

While the H$_2$ concentration of the gas supplied to D1 is 0.1, that of the blend delivered to D2 is in the neighborhood of 0.06 as a result of the mixing of gases from pipes P3 and P4. Moreover, D3 is a fixed demand consumer with no carbon emissions offset, thus the solution does not seek to deliver H$_2$ to at J5.  Moreover, we observe that the shadow prices computed at the two physical nodes differ, which reflects the difference in the locational values of energy, the constituent gases, and carbon emissions mitigation. 

\section{Discussion}  \label{sec:discussion}

We have performed extensive sensitivity analyses on the problem \eqref{prob:hgto} applied to a single pipe in order to demonstrate the intuitive changes to the solution with respect to key parameters.  In analysis \ref{sec:analysis1}, given the pricing structure, no hydrogen injection occurs because hydrogen has lower energy content at the mixture's common pressure, and thus there is no incentive for H$_2$ injection to appear in the solution. This scenario demonstrates the behavior of the system as demand is increased, constraints bind, the pipeline is conjested, and the shadow prices exhibit jumps and ramps to reflect how energy value increases with the demand.

In analysis \ref{sec:analysis2}, the first analysis is repeated but now with an incentive to replace NG with H$_2$ in order to mitigate carbon emissions.   We see that when the system is not congested, i.e., the pipeline is not constrained by pressure and/or compressor constraints at both sending and receiving ends, H$_2$ is supplied to the maximum possible extent.  As energy demand increases, the pipeline becomes congested, and because the value of energy is greater than the value of carbon mitigation, there is a transition to pure NG utilization.

In analysis \ref{sec:analysis3}, the demand is fixed, and the setting reflects an increasing H$_2$ utilization mandate.  Energy transport remains in the form of NG to the extent possible.  Eventually, delivery of energy at the desired level cannot be maintained.  Counter-intuitively, the shadow price of H$_2$ jumps and then gradually decreases between changes in the binding constraint set.  This is evidence of the nonlinearity and non-monotonicity of heterogeneous gas transport. This also indicates that mandating H$_2$ blending could lead to complex market structures with unexpected outcomes.

Finally, the analysis \ref{sec:analysis4} provides an insightful result.  There the energy demand upper bound is given, and the carbon offset price is increased.  We observe increasing shadow prices of both NG and H$_2$ with an increasing carbon offset price, and the solutions show monotone change with respect to this parameter.  When the value of offsetting emissions becomes greater than the value of delivered energy, the gas blend quickly transitions to H$_2$ injection at the maximum possible level subject to engineering constraints.  Future studies could examine this break point in H$_2$ concentration injection to quantify the appropriate ``green premium'' for H$_2$ blending, and examine whether this quantity can be analytically derived from the Karush-Kuhn-Tucker conditions.  Potentially, a locational carbon offset value could be derived, similarly to the notion of locational marginal carbon intensity that was proposed for power grid generation \cite{ruiz2010analysis}.  The analysis of the 8-node case in Section \ref{sec:8node} shows that the valuations are very much location-dependent.  The problem \eqref{prob:hgto} could be used as an optimization-based market mechanism for locational valuation of the emissions reduction achieved with hydrogen injection into natural gas pipelines.  In such a market, suppliers could provide commodity prices of natural gas and hydrogen, and consumers could provide price and quantity bids for energy, as well as the amount they are willing to pay for carbon emissions reduction. 

\section{Conclusions}  \label{sec:conclusion}

We present an economic optimization problem for allocating the flow of natural gas and hydrogen blends through transportation pipeline networks, accounting for delivered energy in withdrawn flows, the cost of natural gas and hydrogen injections, and avoided carbon emissions.  We examine the sensitivity of the physical and dual solution to several parameters to examine carbon mitigation mechanisms. Future work could compare locational and global carbon pricing, and analytically derive dual variables for pricing solutions.

\addtolength{\textheight}{-12cm}   



\vspace{-1ex}
\bibliographystyle{unsrt}  

\typeout{}
\bibliography{refs}

\end{document}